\documentclass[12pt]{article}

\usepackage{amssymb,dsfont}
\usepackage{amsmath}
\usepackage{verbatim}

\usepackage{latexsym}
\usepackage{ifthen}

\usepackage[english]{babel}

\usepackage{a4wide}

\newtheorem{theorem}{Theorem}[section]

\newtheorem{proposition}[theorem]{Proposition}

\newenvironment{Proof}{\vspace{1.5ex}{\sc Proof}.}{\vspace{1.5ex}}

\def\eop{\rule{0.7ex}{0.7ex}}

\newif\ifdisplay\displayfalse

\def\eq{\displaytrue$$}
\def\qe{$$\displayfalse}

\def\eqd{:=}

\def\Re{\mathds{R}}

\def\M{\mathcal{M}}

\def\L{\mathcal{L}}

\def\P{\mathcal{P}}
\def\A{\mathcal{A}}

\def\norm#1{
	\ifdisplay\left\|#1\right\|
	\else\|#1\|\fi
	}

\def\scalf#1#2{
	\ifdisplay\left\langle #1,#2 \right\rangle
	\else\langle #1,#2 \rangle_F\fi
	}
\def\distf{{\mathrm{dist}}_F}

\def\normf#1{
	\ifdisplay\left\|#1\right\|
	\else\|#1\|_F\fi
	}
\def\set#1#2{\{#1|#2\}
    }

\def\tr#1{\mathop{\mathrm{tr}}#1}

\def\argmin{\mathop{\mathrm{argmin}}}

\def\segoo#1#2{(#1, #2)}

\def\vgd{\mathbf{d}}
\def\vgx{\mathbf{x}}

\begin{document}

\title{A proximal approach to the inversion of ill-conditioned matrices}

\date{}

\author{Pierre Mar\'echal\thanks{Institut de Math\'ematiques, Universit\'e Paul Sabatier, 31062 Toulouse, France. {\it pr.marechal@gmail.com}} \and Aude Rondepierre\thanks{Institut de Math\'ematiques, INSA de Toulouse, D\'epartement GMM, 31077 Toulouse, France. {\it aude.rondepierre@math.univ-toulouse.fr}}}
\maketitle

\begin{abstract}
We propose a general proximal algorithm for the inversion of ill-conditioned matrices.
This algorithm is based on a variational characterization of pseudo-inverses.
We show that a particular instance of it (with constant {\sl regularization parameter})
belongs to the class of {\sl fixed point} methods. Convergence of the
algorithm is also discussed.

\vskip 0.5\baselineskip
\end{abstract}

\section{Introduction}
\label{introduction}
Inverting ill-conditioned large matrices is a 
challenging problem involved in a wide range of applications, including
inverse problems (image reconstruction, signal analysis, etc.)
and partial differential equations (computational fluid dynamics, mechanics, etc.).
There are two classes of methods: the first one involves factorization
of the matrix (SVD, QR, LU, LQUP); the second one involves iterative
schemes (fixed point methods, projection onto increasing sequences of subspaces).

The main purpose of this note is to
show that a particular instance of the Proximal Point Algorithm provides
a fixed point method for the problem of matrix inversion.
This fact is based on the observation that the
pseudo-inverse $M^\dagger$ of a matrix $M\in \Re^{m\times n}$
satisfies the {\sl fixed point equation}:
$$
\Phi=\varphi(\Phi)\eqd B\Phi+C,
\quad\hbox{with}\quad
B\eqd (I+\mu M^\top M)^{-1}
\quad\hbox{and}\quad
C\eqd (M^\top M+\mu^{-1}I)^{-1}M^\top
$$
where $\mu>0$. The corresponding {\sl fixed point iteration} $\Phi_{k+1}=B\Phi_k+C$
is nothing but a proximal iteration. We see that $\varphi$ is a
contraction and that, if $M^\top M$ is positive definite,
then $\varphi$ is a strict contraction.
It is worth noticing that, in the proximal algorithm, $\mu$
may depend on~$k$, allowing for large (but inaccurate) steps
for early iteration and small (but accurate) steps when
approaching the solution.

The Proximal Point Algorithm (PPA) was introduced in 1970 by Martinet~\cite{martinet},
in the context of the regularization of variational inequalities.
A few years later, Rockafellar~\cite{rockafellar} generalized the PPA to
the computation of zeros of a maximal monotone operator.
Under suitable assumptions, it can be used to efficiently minimize
a given function, by finding iteratively a zero in its Clarke subdifferential.

Throughout, we denote by $\normf{M}$ the Frobenius norm of a matrix $M\in \Re^{m\times n}$
and by $\scalf{M}{N}$ the Frobenius inner product of $M,N\in\Re^{m\times n}$
(which is given by $\scalf{M}{N}=\tr{(MN^\top)}=\tr{(N^\top M)})$.
In $\Re^{m\times n}$, we denote by $\distf(M,{\mathcal S})$
the distance between a matrix~$M$ and a set ${\mathcal S}$:
$$
\distf(M,{\mathcal S})\eqd\inf\set{\normf{M-M'}}{M'\in {\mathcal S}}.
$$
The identity matrix will be denoted by~$I$, its dimension being always clear from 
the context.

The next theorem, whose proof may be found {\it e.g.} in~\cite{amodei-dedieu},
provides a variational characterization of~$M^\dagger$.

\begin{theorem}
\label{theorem1}
The pseudo-inverse of a matrix $M\in \Re^{m\times n}$
is the solution of minimum Frobenius norm of the optimization
problem
\eq
(\P)\quad
\hbox{Minimize}\;
f(\Phi)\eqd\frac{1}{2}\norm{M\Phi-I}_F^2
\;\;\hbox{over}\;\;
\Re^{n\times m}.
\qe
\end{theorem}

\section{The proximal point algorithm}
\label{proximal}

The proximal point algorithm is a general algorithm for computing zeros
of maximal monotone operators. A well-known application is the minimization
of a convex function $f$ by finding a zero in its subdifferential.
In our setting, it consists in the following steps:
\begin{enumerate}
\item[1.] Choose an initial matrix $\Phi_0\in \Re^{m\times n}$;
\item[2.] Generate a sequence $(\Phi_k)_{k\geq 0}$ according to the formula
\begin{equation}
\label{proximal:iteration}
\Phi_{k+1}=\displaystyle\argmin_{\Phi\in \Re^{m\times n}}\left\{f(\Phi)
+\frac{1}{2\mu_k}\normf{\Phi-\Phi_k}^2\right\},
\end{equation}
in which $(\mu_k)_{k\geq 0}$ is a sequence of positive numbers,
until some stopping criterion is satisfied.
\end{enumerate}

Equation~\eqref{proximal:iteration} will be subsequently referred to as the
{\sl proximal iteration} of Problem $(\P)$. The stopping criterion may combine, as usual, conditions such as 
$$
\normf{\nabla f(\Phi_k)}\leq\varepsilon_1
\quad\hbox{and}\quad
\normf{\Phi_k-\Phi_{k-1}}\leq\varepsilon_2,
$$
where the parameters $\varepsilon_1$ and $\varepsilon_2$ control the precision
of the algorithm.

Clearly, the function $f:\Phi\mapsto\normf{M\Phi-I}^2/2$ is convex and
indefinitely differentiable. Therefore, solutions of the proximal iteration
(\ref{proximal:iteration}) are characterized by the relationship
$\nabla f(\Phi_{k+1})+\mu_k^{-1}(\Phi_{k+1}-\Phi_k)=0$ i.e.,
\begin{equation}
\label{prox2}
(I + \mu_k M^\top M)\Phi_{k+1} = \Phi_k + \mu_k M^\top.
\end{equation}
Since $M^\top M$ is positive semi-definite and $\mu_k$ is chosen to be positive
for all $k$, the
matrix $(I+\mu_kM^\top M)$ is nonsingular and the proximal iteration also reads:
\begin{equation}
\label{prox3}
\Phi_{k+1}=\left(I+\mu_k M^\top M\right)^{-1}\left(\Phi_k+\mu_k M^\top\right).
\end{equation}
The iterates $\Phi_k$ could be computed either exactly (in the ideal case),
or approximately, using {\it e.g.} any efficient minimization algorithm to solve
the proximal iteration~\eqref{proximal:iteration}. In that case, we need another
stopping criterion and we here choose the following one suggested in \cite{luque}:
\begin{equation}
\normf{\Phi_{k+1}-A(\Phi_k +\mu_k M^\top)} \leq
\epsilon_k \min\{1,\normf{\Phi_{k+1}-\Phi_k}^r\},~r>1\label{Ar}
\end{equation}
where $\epsilon_k>0$ and the series $\sum\epsilon_k$ is convergent.
Notice that the larger $r$, the more accurate the computation of $\Phi_{k+1}$.
Notice also that, in the case where $\mu_k=\mu$ for all $k$, each proximal
iteration involves the multiplication by the same matrix $A\eqd(I+\mu M^\top M)^{-1}$,
and that the latter inverse may be easy to compute numerically, if the
matrix $I+\mu M^\top M$ is well-conditioned.

We now turn to convergence issues.
Recall that our objective function~$f$ is a quadratic function
whose Hessian $M^\top M$ is positive semi-definite. Nevertheless,
unless $M^\top M$ is positive definite, the matrix $I-(I+\mu M^\top M)^{-1}$
is in general singular and the classical convergence theorem for iterative
methods (see e.g. \cite{Ciarlet}) is not helpful here to prove the convergence
of our proximal scheme.
The following proposition is a consequence of Theorem~2.1 in~\cite{luque}.
For clarity, we shall denote by $\M$ the linear mapping $\Phi\mapsto M\Phi$,
by~$\L$ the linear mapping $\Phi\mapsto M^\top M\Phi$ and
by $\A$ the linear mapping $\Phi\mapsto A\Phi=(I+\mu M^\top M)^{-1}\Phi$.

\begin{proposition}
\label{convergence-result}
Let $\alpha_1$ be the smallest nonzero eigenvalue of $\L$ and let $E_1$
be the corresponding eigen\-space.
Assume that $\mu_k=\mu$ for all $k$ and that $\Phi_0$ is not $\scalf{\cdot}{\cdot}$-orthogonal
to the eigenspace $E_1$. Then,
$$
\frac{\normf{A(\Phi_{k+1}-\Phi_k)}}{\normf{\Phi_{k+1}-\Phi_k}}
\to\frac{1}{1+\alpha_1\mu}
\quad\hbox{and}\quad
\frac{\Phi_{k+1}-\Phi_k}{\normf{\Phi_{k+1}-\Phi_k}}
\to\Psi_1
\quad\hbox{as}\quad k\to\infty,
$$
in which $\Psi_1$ is a unit eigenvector in $E_1$.
Moreover the sequence $(\Phi_k)$ generated by the proximal algorithm,
either with infinite precision or using the stopping criterion~\eqref{Ar}
for the inner loop, converges linearly to the orthogonal projection of~$\Phi_0$
onto the solution set ${\mathcal S}\eqd\argmin{f}=M^\dagger+\ker\M$.
\end{proposition}

\Proof{
{\sl Step 1}. The error $\Delta_{k+1}\eqd\Phi_{k+1}-\Phi_k$ (at iterate $k+1$)
satisfies: $\Delta_{k+1}=(I+\mu M^\top M)^{-1}\Delta_k$.
The latter iteration is that of the power method for the linear mapping~$\A$.
Clearly, ${\mathcal A}$ is symmetric and positive definite. Consequently,
its eigenvalues are strictly positive and there exists an unique eigenvalue
of largest modulus (not necessarily simple).

Notice now that the eigenspace $E_1$ associated to the eigenvalue $\alpha_1$
is nothing but the eigenspace of ${\mathcal A}$ for its largest eigenvalue
strictly smaller than~1, namely, $1/(1+\alpha_1\mu)$. We proceed as in~\cite[Theorem 1]{Vige}
to obtain the desired convergence rate {\it via} that of the iterated power method. 

{\sl Step 2}. We now establish the linear convergence of the sequence $(\Phi_k)$.
First, the solution set~${\mathcal S}$ is clearly nonempty since it contains~$M^\dagger$.
Moreover, let $(\Phi_k)$ be a sequence generated by the PPA algorithm using the
stopping criterion \eqref{Ar}. Let us prove that
\begin{equation}
\exists a>0,\;
\exists \delta>0,\;
\forall\Phi\in \Re^{m\times n},\quad
\Bigl[\normf{\nabla f(\Phi)}<\delta\Rightarrow
\distf(\Phi,{\mathcal S})\leq a\normf{\nabla f(\Phi)}\Bigr]\label{hyp:cv},
\end{equation}
which is nothing but Condition (2.1) in \cite[Theorem 2.1]{luque},
in our context.
Let $\Phi\in \Re^{m\times n}$ and let $\bar{\Phi}$ be the orthogonal
projection of $\phi$ over $(\ker\M)^\perp$.
It results from the classical theory of linear least squares that
$\distf(\Phi,{\mathcal S})=\normf{\bar{\Phi}-M^\dagger}$. 
Since $M^\top M \bar\Phi-M^\top=\nabla f(\bar\Phi)$
and $M^\top M M^\dagger-M^\top=0$, we also have:
$\nabla f(\Phi)=M^\top M(\bar\Phi-M^\dagger)$. 
Moreover, $\bar\Phi-M^\dagger\in(\ker\M)^\perp=(\ker\L)^\perp$,
so that
$$
\normf{\nabla f(\Phi)}=
\normf{M^\top M(\bar\Phi-M^\dagger)}\geq \alpha_1 \normf{\bar\Phi-M^\dagger}.
$$
It follows that~\eqref{hyp:cv} is satisfied with $a = 1/\alpha_1$.
The conclusion then follows from \cite[Theorem 2.1]{luque}:
the sequence $(\Phi_k)$ converges linearly with a rate bounded by
$a/\sqrt{a^2+\mu^2}=1/\sqrt{1+\mu^2\alpha_1^2}<1$.

{\sl Step 3}. By rewriting the proximal iteration in an orthonormal
basis of eigenvectors of~${\mathcal L}$, we finally prove that the
limit of the sequence $(\Phi_k)$ is the orthogonal projection
of $\Phi_0$ onto $\argmin{f}$.~\eop}

A complete numerical study, which goes beyond the scope of this paper,
is currently in progress and will be presented in a forthcoming
publication. Let us merely mention that our proximal approach makes it
possible to combine features of factorization methods (in the proximal
iteration) with features of iterative schemes. In particular, if~$M$
if invertible, it shares with iterative methods the absence of error
propagation and amplification, since each iterate can be regarded as
a new initial point of a sequence which converges to the desired solution.

\section{Comments}
\noindent
{\bf Tikhonov approximation.}
A standard approximation of the pseudo-inverse of an ill-conditioned
matrix~$M$ is $(M^\top M+\varepsilon I)^{-1}M^\top$, where $\varepsilon$
is a small positive number. This approximation
is nothing but the Tikhonov regularization of~$M^\dagger$, with regularization
parameter~$\varepsilon$. It is worth noticing
that the choice $\Phi_0=0$ in the proximal algorithm yields the latter approximation
for $\varepsilon= 1/\mu$ after one proximal iteration.
\vspace{1ex}

\noindent
{\bf Trade-offs.}
At the $k$-th proximal iteration, the perturbation of the objective
function~$f$ is, roughly speaking, proportional to the square of the
distance between the current iterate and the solution set of~$(\P)$,
and inversely proportional to $\mu_k$.
In order to speed up the algorithm, it seems reasonable to choose
large $\mu_k$ for early iterations, yielding large but inaccurate
steps, and then smaller $\mu_k$ for late iterations, where {\sl proximity}
with the solution set makes it suitable to perform small and accurate
steps.
This is especially true in the case where $M$ is invertible, since
the solution set then reduces to~$\{M^{-1}\}$.
Moreover, numerical accuracy in early proximal iteration may be
irrelevant, since the limit of the proximal sequence is what
really matters. A trade-off between a rough approximation of the
searched proximal point and an accurate and costly solution
must be found. As suggested in \cite{CL:93}, one may use the following
stopping criterion for the proximal iteration:
$$
f(\Phi_{k+1})-f(\Phi_k)\leq\delta\scalf{\nabla f(\Phi_{k+1})}{\Phi_{k+1}-\Phi_k}.
$$
This criterion is an Armijo-like rule: the algorithm stops when
the improvement of the objective function~$f$ is at least a given
fraction $\delta\in\segoo{0}{1}$ of its ideal improvement.
\vspace{1ex}


\noindent
{\bf Inversion versus linear systems.}
It is often unnecessary to compute the inverse of a matrix~$M$, in particular
when the linear system $M\vgx=\vgd$ must be solved for a few data vectors $\vgd$ only.
In such cases, of course, the usual proximal strategy may be used to compute least
squares solutions. It is important to realize that, although the regularization
properties of the proximal algorithm are effective at every proximal iteration,
perturbations of $\vgd$ may still have dramatic effects on the algorithm
if $M$ is ill-conditioned. In applications for which no perturbation of the data
must be considered, accurate solutions may be reached by a proximal strategy.
We emphasize that, in the minimization of $\Phi\mapsto\normf{M\Phi-I}$, the  
{\sl data}~$I$ undergoes no perturbation whatsoever.
\vspace{1ex}



\begin{thebibliography}{10}

\bibitem{amodei-dedieu}
L. Amodei and J.-P. Dedieu,
\newblock {\em Analyse Num\'erique Matricielle}.
\newblock Collection Sciences Sup, Dunod, 2008.

\bibitem{Ciarlet} P.G. Ciarlet.
\newblock {\em Introduction to numerical linear algebra and optimisation}.
\newblock Cambridge Texts in Applied Mathematics, Cambridge University Press, 1989.

\bibitem{CL:93} R. Correa and C. Lemar{\'e}chal.
\newblock {\em Convergence of some algorithms for convex minimization}.
\newblock Math. Programming, vol. 62, pp. 161--275, 1993.

\bibitem{luque} F.J. Luque.
\newblock Asymptotic convergence analysis of the proximal point algorithm.
\newblock SIAM Journal on Control and Optimization, vol. 22 (2), pp. 277--293, 1984.

\bibitem{martinet}
B. Martinet.
\newblock{\em R\'egularisation d'in\'equations variationelles par approximations successives}.
\newblock Revue Fran\c caise d'Informatique et de Recherche Op\'erationelle,
pp. 154--159, 1970.

\bibitem{rockafellar}
R.T. Rockafellar.
\newblock{\em Monotone operators and the proximal point algorithm}.
\newblock SIAM Journal on Control and Optimization, vol. 14 (5), p877--898, 1976. 

\bibitem{Vige} G. Vige
\newblock {\em Proximal-Point Algorithm for Minimizing Quadratic Functions}.
\newblock INRIA research report RR-2610, 1995.
\end{thebibliography}
\end{document}